\documentclass[11pt]{article}

\usepackage{latexsym,amsmath,amsthm,amsfonts,amsbsy,amscd}
\theoremstyle{plain}
\newtheorem{lemma}{Lemma}[section]

\newtheorem{proposition}{Proposition}[section]
\newtheorem{corollary}{Corollary}[section]
\newtheorem{definition}{Definition}[section]

\numberwithin{equation}{section}

\newcommand{\qbin}[2]{\genfrac{[}{]}{0pt}{}{#1}{#2}}
\newcommand{\qbins}[2]{{\textstyle\genfrac{[}{]}{0pt}{}{#1}{#2}}}
\def\Integer{\mathbb{Z}}
\def\I{{\cal I}}
\def\ng{\negthickspace}

\begin{document}

\title{$q$-Trinomial identities}
 
\author{\Large
S.~Ole Warnaar\thanks{
Instituut voor Theoretische Fysica,
Universiteit van Amsterdam,
Valckenierstraat 65,
1018 XE Amsterdam, The Netherlands;
e-mail: {\tt warnaar@phys.uva.nl}}}
 
\maketitle

\begin{abstract}
We obtain connection coefficients between $q$-binomial and
$q$-trinomial coefficients. Using these, one can transform
$q$-binomial identities into a $q$-trinomial identities and
back again.
To demonstrate the usefulness of this procedure we rederive some
known trinomial identities related to partition theory
and prove many of the conjectures of Berkovich, McCoy and Pearce,
which have recently arisen in their study of the $\phi_{2,1}$ and $\phi_{1,5}$
perturbations of minimal conformal field theory. 
\end{abstract}

\section{Introduction}
The $q$-binomial coefficients can be defined by the $q$-analogue of
Newton's binomial expansion,
\begin{equation}\label{qNewton}
(1+x)(1+qx)\dots(1+q^{n-1}x)=\sum_{a=0}^n x^a q^{a(a-1)/2}\qbin{n}{a}.
\end{equation}
An explicit expression for the $q$-binomial coefficients is given by
\begin{equation*}
\qbin{n}{a}_q=\qbin{n}{a}=
\begin{cases}\displaystyle \frac{(q)_{n}}{(q)_a(q)_{n-a}} & 
\text{for $0\leq a \leq n$} \\[3mm]
0 & \text{otherwise,} \end{cases}
\end{equation*}
where 
\begin{equation*}
(q)_n=\prod_{j=1}^n (1-q^j)\quad n\geq 1 \qquad \text{and} \qquad
(q)_0=1.
\end{equation*}

$q$-Binomials play an essential role in combinatorics,
partition theory and statistical mechanics, see e.g., 
\cite{Andrews76,Andrews85,GR90,MacMahon16}, and one of 
MacMahon's famous results is that $\qbin{n+m}{m}$ is the
generating function of partitions with no more than $m$ parts, no
part exceeding $n$.
Less well-understood are the $q$-trinomial coefficients, defined as 
$q$-analogues of the numbers appearing in the generalized Pascal triangle
\begin{equation}\label{Pascal}
\begin{array}{ccccccccc}
&&&&1&&&& \\
&&&1&1&1&&& \\
&&1&2&3&2&1&& \\
&1&3&6&7&6&3&1& \\
.&.&.&.&.&.&.&.&.
\end{array}
\end{equation}
Andrews and Baxter were the first to encounter $q$-trinomial coefficients,
and in ref.~\cite{AB87} they defined
\begin{equation}\label{qt}
\qbin{L,b;q}{a}_2=\qbin{L,b}{a}_2
=\sum_{k\geq 0} q^{k(k+b)}\qbin{L}{k}\qbin{L-k}{k+a}
\end{equation}
and
\begin{equation}\label{qtT}
T_n(L,a;q)=T_n(L,a)=
q^{\frac{1}{2}(L-a)(L+a-n)}\qbin{L,a-n;q^{-1}}{a}_2.
\end{equation}
The $q$-trinomial $T_n$ can be expressed explicitly as 
\begin{equation}\label{qtTn}
T_n(L,a)=\sum_{\substack{r=0 \\L-a-r \text{ even}}}^{L-|a|}
\frac{q^{\frac{1}{2}r(r-n)}(q)_L}{(q)_{\frac{L-a-r}{2}}
(q)_{\frac{L+a-r}{2}}(q)_r}\, .
\end{equation}
Clearly, the $q$-trinomials coefficients are nonzero for 
$a=-L,-L+1,\dots,L$ only and satisfy the symmetries
\begin{equation*}
\qbin{L,b;q}{a}_2=q^{a(a-b)}\qbin{L,b-2a}{-a}_2 \qquad \text{and}
\qquad T_n(L,a)=T_n(L,-a).
\end{equation*}
To see that \eqref{qt} indeed defines $q$-analogues of the trinomial
coefficients, set $q=1$ and twice apply the binomial formula to find that
\begin{equation*}
\sum_{a=-L}^L x^a \qbin{L,b;1}{a}_2=(1+x+x^{-1})^L 
\end{equation*}
in accordance with \eqref{Pascal}.
The only further properties of $q$-trinomials needed in this paper
are the limiting formulas~\cite{AB87}
\begin{align}\label{Tlimeven}
\lim_{\substack{L\to\infty \\L-a\text{ even}}} T_0(L,a)&=
\frac{(-q^{1/2})_{\infty}+(q^{1/2})_{\infty}}{2(q)_{\infty}} \\
\label{Tlimodd}
\lim_{\substack{L\to\infty \\L-a\text{ odd}}} T_0(L,a)&=
\frac{(-q^{1/2})_{\infty}-(q^{1/2})_{\infty}}{2(q)_{\infty}}
\end{align}
and
\begin{equation}\label{trinlim}
\lim_{L\to\infty} \qbin{L,a}{a}_2=\frac{1}{(q)_{\infty}}.
\end{equation}
Finally we introduce the abbreviation
\begin{equation*}
\qbin{L,a}{a}_2=\qbin{L}{a}_2.
\end{equation*}

Since their discovery about a decade ago, $q$-trinomials have found
numerous applications in, again, combinatorics, partition theory,
and statistical mechanics~[4--8,11,16--19,32,33,35,38--42].
Among the most striking results is a $q$-trinomial proof of Schur's 
partition theorem and Capparelli's (then) conjecture~\cite{Andrews94},
a $q$-trinomial proof of the G\"ollnitz--Gordon partition 
theorem~\cite{Andrews90a} and their Andrews--Bressoud 
generalization's~\cite{BM97,BMO96}, the proof of an E$_8$ 
Rogers--Ramanujan-type identity~\cite{WP94}, and a trinomial analogue of
Bailey's lemma~\cite{AB98}.

Most of the above-cited papers contain $q$-trinomial identities.
Upon close inspection of many of these identities
one is struck by their similarity with well-known $q$-binomial identities.
This strongly suggests that many $q$-trinomial identities can be simply
viewed as corollaries of $q$-binomial identities.
In an earlier paper~\cite{Warnaar98} we made a first, only partially 
successful, attempt to relate $q$-trinomial identities to $q$-binomial 
identities, showing that each Bailey pair (which implies a $q$-binomial 
identity) implies a trinomial Bailey pair (which implies a $q$-trinomial 
identity).  The problem with the idea of \cite{Warnaar98} is that it applies
to $q$-trinomial identities in which the parameter $a$ in \eqref{qt} and 
\eqref{qtT} takes even values only. Therefore, $q$-trinomial identities in 
which $a$ takes arbitrary integer values remained irreducible to 
$q$-binomial identities.

This papers is intended to deal with this problem, and in the next
section we obtain connection coefficients between $q$-binomial
and $q$-trinomial coefficients.
Using these coefficients and the idea of ref.~\cite{Warnaar98}
many $q$-trinomial identities are derived from known $q$-binomial
identities. In section~\ref{secpt} several $q$-trinomial identities
related to partitions are obtained and in section~\ref{sectid}
we prove general classes of $q$-trinomial identities,
including many of the recent conjectures of Berkovich, 
McCoy and Pearce~\cite{BMP98}.
To make contact with the recently discovered trinomial analogue of Bailey's
lemma, we finally formulate our results in the language of Bailey pairs 
in section~\ref{secBl}. In the appendix some necessary formulas
for $q$-binomial coefficients are collected.

\section{Connection coefficients}\label{seccon}
To relate $q$-binomials and $q$-trinomials we consider the simple problem 
of finding the coefficients $C_{L,k}$ and $C'_{L,k}$ such that
\begin{equation}\label{TB}
T_0(L,a)=\sum_{k=0}^L C_{L,k}(a)\qbin{2k}{k-a}
\end{equation}
and
\begin{equation}\label{BT}
\qbin{2L}{L-a}=\sum_{k=0}^L C'_{L,k}(a)T_0(k,a).
\end{equation}
Of course, the two equations imply that
\begin{equation}\label{ortho}
\sum_{k=M}^L C_{L,k}(a) C'_{k,M}(a)=\delta_{L,M}.
\end{equation}

The answer to the above connection coefficient problem is given by
the following lemma.
\begin{lemma}\label{lemma1}
For $C_{L,k}$ and $C'_{L,k}$ as above,
\begin{align}\label{C}
C_{L,k}(a)&=(-1)^{L-k}q^{\binom{L-k}{2}+\frac{1}{2}(a^2-L^2)}\qbin{L}{k} \\
\label{Cp}
C'_{L,k}(a)&=q^{\frac{1}{2}(k^2-a^2)}\qbin{L}{k}.
\end{align}
\end{lemma}
\begin{proof}
Substitution of the expression for $C'_{L,k}$ into the right-hand side of
\eqref{BT} and using equation \eqref{qtTn} for $T_0$ gives
\begin{equation*}
\sum_{k=0}^L C'_{L,k}(a)T_0(k,a)=
\sum_{k=0}^L \sum_{\substack{r=0 \\k-a-r \text{ even}}}^{k-|a|}
\frac{q^{\frac{1}{2}(k^2-a^2+r^2)}(q)_L}{(q)_{L-k}(q)_{\frac{k-a-r}{2}}
(q)_{\frac{k+a-r}{2}}(q)_r}.
\end{equation*}
To proceed we introduce new summation variables $i,j$ defined by $k=i+j+a$ 
and $r=i-j$, and apply the $q$-Chu--Vandermonde sum, i.e.,
\begin{eqnarray*}
\sum_{k=0}^L C'_{L,k}(a)T_0(k,a) \ng\ng &=& \ng\ng
\sum_{i=0}^L \sum_{j=0}^i q^{i(i+a)+j(j+a)}
\qbin{L}{i}\qbin{i}{j}\qbin{L-i}{j+a} \notag \\
\ng\ng &\stackrel{\text{by }\eqref{qcv1}}{=}& \ng\ng
\sum_{i=0}^L q^{i(i+a)}\qbin{L}{i}\qbin{L}{i+a} \notag \\[1.5mm]
\ng\ng 
&\stackrel{\text{by }\eqref{qcv1}}{=}& \ng\ng \qbin{2L}{L-a}.
\end{eqnarray*}
This settles \eqref{Cp} and to prove \eqref{C} we show that \eqref{ortho} 
holds. Taking the left-hand side of \eqref{ortho} and substituting the 
claim of the lemma we find
\begin{align*}
\sum_{k=M}^L C_{L,k}(a) C'_{k,M}(a)&=
\sum_{k=M}^L (-1)^{L-k} q^{\binom{L-k}{2}+\frac{1}{2}(M^2-L^2)}\qbin{L}{k} 
\qbin{k}{M} \notag \\
&=q^{\frac{1}{2}(M^2-L^2)}\qbin{L}{M}
\sum_{k=0}^{L-M}(-1)^k q^{\binom{k}{2}}\qbin{L-M}{k} =\delta_{L,M},
\end{align*}
where in the last step we have used \eqref{qNewton} with $x=-1$.
\end{proof}

The analogous result involving $T_1$ instead of $T_0$ can be stated as follows.
Define $D_{L,k}$ and $D'_{L,k}$ by
\begin{equation*}
T_1(L,a)=\sum_{k=0}^L D_{L,k}(a)\qbin{2k}{k-a}
\end{equation*}
and
\begin{equation}\label{BT1}
\qbin{2L}{L-a}=\sum_{k=0}^L D'_{L,k}(a)T_1(k,a).
\end{equation}
\begin{lemma}\label{lemma2}
For $D_{L,k}$ and $D'_{L,k}$ as above,
\begin{align}\label{D}
D_{L,k}(a)&=
(-1)^{L-k}q^{\binom{L-k}{2}+\binom{a}{2}-\binom{L}{2}}
\frac{1+q^a}{1+q^k}\qbin{L}{k} \\
\label{Dp}
D'_{L,k}(a)&=q^{\binom{k}{2}-\binom{a}{2}}\frac{1+q^L}{1+q^a}\qbin{L}{k}.
\end{align}
\end{lemma}
\begin{proof}
Following the proof of lemma~\ref{lemma1} with $T_0$ replaced by $T_1$, one 
finds after application of the $q$-Chu--Vandermonde sum \eqref{qcv1}, that 
the right-hand side of \eqref{BT1} is equal to
\begin{equation*}
\frac{1+q^L}{1+q^a}
\sum_{i=0}^L q^{i(i+a-1)}\qbin{L}{i}\qbin{L}{i+a}.
\end{equation*}
Before \eqref{qcv1} can again be applied the recurrence \eqref{qbinrec} is
needed to rewrite this as
\begin{equation*}
\frac{1+q^L}{1+q^a} \Bigl\{
\sum_{i=0}^L q^{i(i+a-1)}\qbin{L}{i}\qbin{L-1}{i+a-1}
+q^a\sum_{i=0}^L q^{i(i+a)}\qbin{L}{i}\qbin{L-1}{i+a} \Bigr\}.
\end{equation*}
Using \eqref{qcv1} and combining terms gives $\qbins{2L}{L-a}$.
To prove \eqref{D} it again suffices to consider
$\sum_{k=M}^L D_{L,k}(a) D'_{k,M}(a)$. After substituting the results for 
$D$ and $D'$ and replacing $k\to L-k$ one finds that this becomes 
$\delta_{L,M}$ after using \eqref{qNewton} with $x=-1$.
\end{proof}

To conclude this section we note that the representations \eqref{qt} and
\eqref{qtTn} for the $q$-trinomial coefficients
can also be written as a relation between $q$-trinomials and $q$-binomials.
That is, 
\begin{equation}\label{TBe}
T_n(L,2a)=\sum_{k\geq 0}
q^{\frac{1}{2}(L-2k)(L-2k-n)} \qbin{L}{2k}\qbin{2k}{k-a}
\end{equation}
and
\begin{equation}\label{tBe}
\qbin{L,b}{2a}_2=\sum_{k\geq 0}
q^{(k-a)(k-a+b)}\qbin{L}{2k}\qbin{2k}{k-a}.
\end{equation}
These results, which unlike the previous transformations are not invertible,
will be needed later.

\section{Simple examples from partition theory}\label{secpt}
Before proving general series of $q$-trinomial identities using the 
results of the previous section, we treat some simple examples related 
to partition identities first.

The first example concerns the following result of 
Andrews~\cite{Andrews90a} (see also~\cite{BMP98}).
It is well known~\cite{Andrews70} that the (first) Rogers--Ramanujan 
identity can be obtained as limiting case of the polynomial identity
\begin{equation}\label{RR}
\sum_{n \geq 0}q^{n^2}\qbin{L-n}{n}=
\sum_{j=-\infty}^{\infty}(-1)^j q^{j(5j+1)/2}
\qbin{L}{\lfloor \frac{L-5j}{2} \rfloor}.
\end{equation}
Here the polynomials appearing on either side are known to be the 
generating function of partitions with difference between parts at least two 
and largest part not exceeding $L-1$~\cite{MacMahon16,Schur17}.
In ref.~\cite{Andrews90a}, Andrews remarks that it is ``most 
surprising and intriguing'' that the following also holds
\begin{equation}\label{RR2}
\sum_{n\geq 0}q^{n^2}\qbin{L-n}{n}=
\sum_{j=-\infty}^{\infty} \left\{
q^{j(10j+1)}\qbin{L}{5j}_2-
q^{(2j+1)(5j+2)}\qbin{L}{5j+2}_2\right\}.
\end{equation}
We now show that \eqref{RR2} is a corollary of \eqref{RR}, or for those
who prefer to decrease instead of increase complexity, that
\eqref{RR} is a corollary of \eqref{RR2}.
Replacing $q\to 1/q$ in \eqref{TB} and \eqref{C}, using \eqref{qtT} and
\eqref{qbindual}, we find that
\begin{equation*}
\qbin{L}{a}_2=\sum_{k=0}^L (-1)^{L-k}q^{\frac{1}{2}(L-k)(L+k+1)}
\qbin{L}{k}\qbin{2k}{k-a}.
\end{equation*}
Therefore, if we take \eqref{RR} with $L$ replaced by $2k$, multiply 
by $(-1)^{L-k}q^{\frac{1}{2}(L-k)(L+k+1)}\qbins{L}{k}$ and sum over $k$, 
we arrive at
\begin{multline*}
\sum_{k\geq 0}\sum_{n\geq 0}
(-1)^{L-k}q^{\frac{1}{2}(L-k)(L+k+1)+n^2}
\qbin{L}{k}\qbin{2k-n}{n} \\
=\sum_{j=-\infty}^{\infty} \left\{
q^{j(10j+1)}\qbin{L}{5j}_2-
q^{(2j+1)(5j+1)}\qbin{L}{5j+2}_2\right\}.
\end{multline*}
To simplify the left-hand side we set $k=L-m+n$ followed by $n\to m-n$ to get
\begin{equation*}
\sum_{m\geq 0}q^{m^2}
\sum_{n\geq 0} (-1)^n q^{\binom{n}{2}+n(L-2m+1)}
\qbin{L}{n}\qbin{2L-m-n}{m-n}=
\sum_{m\geq 0}q^{m^2}\qbin{L-m}{m},
\end{equation*}
where the sum over $n$ has been performed using
the $q$-Chu--Vandermonde summation~\eqref{qcv3}.
As remarked before, one can equally well take the reverse route and
starting from \eqref{RR2}, using lemma~\ref{lemma1},
one readily obtains \eqref{RR}. We leave this to the reader.

Our second example concerns the following identity of Slater~\cite{Slater52}
related to the (first) G\"ollnitz--Gordon partition 
identity~\cite{Goellnitz67,Gordon65},
\begin{equation}\label{GG}
\sum_{n=0}^{\infty} \frac{q^{n^2}(-q;q^2)_n}{(q^2;q^2)_n}=
\prod_{n=0}^{\infty}\frac{1}{(1-q^{8j+1})(1-q^{8j+4})(1-q^{8j+7})}.
\end{equation}
A polynomial identity that implies this equation is given 
by~\cite{BMO96,BM97}
\begin{equation}\label{BMO}
\sum_{m,n \geq 0}q^{\frac{1}{2}(m^2+n^2)}\qbin{L-m}{n}\qbin{n}{m}
=\sum_{j=-\infty}^{\infty}(-1)^j q^{2j^2+\frac{1}{2}j}\Bigl\{
T_0(L,4j)+T_0(L,4j+1)\Bigr\}. 
\end{equation}
It was observed in \cite{Andrews90a} that for fixed $L$
the polynomial appearing on the right-hand side with $q$ replaced by $q^2$ 
is the generating function of partitions $\lambda=(\lambda_1,\lambda_2,\dots)$
with $\lambda_i-\lambda_{i+1}\geq 2$ for $\lambda_i$ odd, 
$\lambda_i-\lambda_{i+1}\geq 3$ for $\lambda_i$ even, and with largest
part not exceeding $2L-1$.
To see that \eqref{BMO} indeed implies \eqref{GG}, let $L$ tend to infinity
using \eqref{Tlimeven}, \eqref{Tlimodd} and \eqref{qbinlim}. Hence,
\begin{equation*}
\sum_{m,n \geq 0} \frac{q^{\frac{1}{2}(m^2+n^2)}}{(q)_n}\qbin{n}{m}
=\frac{(-q^{1/2})_{\infty}}{(q)_{\infty}}
\sum_{j=-\infty}^{\infty}(-1)^j q^{2j^2+\frac{1}{2}j}.
\end{equation*}
Using Jacobi's triple product identity (equation (2.2.10) of \cite{Andrews76})
and equation \eqref{qNewton} with $x=q^{1/2}$ gives
\begin{equation*}
\sum_{n \geq 0} \frac{q^{\frac{1}{2}n^2}(-q^{1/2})_n}{(q)_n}
=\frac{(-q^{1/2})_{\infty}(q^{3/2};q^4)_{\infty}
(q^{5/2};q^4)_{\infty}(q^4;q^4)_{\infty}}{(q)_{\infty}}.
\end{equation*}
Letting $q\to q^2$ and cleaning up the right-hand side finally yields
\eqref{GG}.

The companion $q$-binomial identity of \eqref{BMO} is given by the
following identity of refs.~\cite{BM96,BMS98}
\begin{multline*}
\sum_{\substack{m_1,m_2 \geq 0 \\m_1+m_2 \text{ even}}}
q^{\frac{1}{4}(m_1^2+m_2^2)}
\qbin{L+\frac{1}{2}(m_1-m_2)}{m_1}\qbin{\frac{1}{2}(m_1+m_2)}{m_2}\\
=\sum_{j=-\infty}^{\infty}(-1)^j \Bigl\{
q^{\frac{1}{2}j(20j+1)}\qbin{2L}{L-4j}+
q^{\frac{1}{2}(4j+1)(5j+1)}\qbin{2L}{L-4j-1}\Bigr\}.
\end{multline*}
To prove this we replace $L$ by $k$, multiply by
$q^{-a^2/2}C_{L,k}(a)$ as given by \eqref{C} and sum over $k$
using \eqref{TB}. The resulting equation is
\begin{eqnarray*}
\lefteqn{
\sum_{j=-\infty}^{\infty}(-1)^j q^{2j^2+\frac{1}{2}j}\Bigl\{
T_0(L,4j)+T_0(L,4j+1)\Bigr\}} \notag \\
\ng\ng&=&\ng\ng\sum_{\substack{m_1,m_2 \geq 0 \\m_1+m_2 \text{ even}}}
q^{\frac{1}{4}(m_1^2+m_2^2-2L^2)}
\qbins{\frac{1}{2}(m_1+m_2)}{m_2}\sum_{k=0}^L(-1)^k 
q^{\binom{k}{2}}\qbins{L}{k}\qbins{L-k+\frac{1}{2}(m_1-m_2)}{m_1}
\notag \\[1.5mm]
\ng\ng&\stackrel{\text{by }\eqref{qcv2}}{=}&\ng\ng
\sum_{\substack{m_1,m_2\geq 0\\m_1+m_2 \text{ even}}}
q^{\frac{1}{4}((m_1-L)^2+(m_2-L)^2)}
\qbins{\frac{1}{2}(m_1+m_2)}{m_2} \qbins{\frac{1}{2}(m_1-m_2)}{m_1-L}
\end{eqnarray*}
Making the variable change $m_1 \to L+n-m$ and $m_2 \to L-n-m$
we find equation \eqref{BMO}.

\section{$q$-Trinomial identities}\label{sectid}
After the previous examples we now derive general
classes of $q$-trinomial identities,
as stated in propositions \ref{prop1}--\ref{prop5} below.
The setup will be as follows. First we describe a family of
$q$-binomial identities for bounded analogues of Virasoro characters,
based on continued fraction expansions.
We then transform these identities into $q$-trinomial identities,
by either using \eqref{TB} or \eqref{tBe}.
Many of the $q$-trinomial identities available in the literature are contained
in the propositions \ref{prop1}--\ref{prop5} or can be derived in
completely analogous fashion.

\subsection{$q$-Binomial identities for bounded Virasoro 
characters}\label{secBM}
Using the inclusion-exclusion construction of Feigin and Fuchs~\cite{FF83},
the (normalized) characters of the Virasoro algebra of
central charge $c=1-6(p'-p)^2/pp'$, with $p,p'$ integers such that 
$1<p<p'$ and $\gcd(p,p')=1$, are given 
by~\cite{RochaCaridi85,Dobrev87}
\begin{equation*}
\chi_{r,s}^{(p,p')}(q)=
\frac{1}{(q)_{\infty}}\sum_{j=-\infty}^{\infty}\Bigl\{
q^{j(pp'j+p'r-ps)}-q^{(pj+r)(p'j+s)}\Bigr\}.
\end{equation*}
Here $r=1,\dots,p-1$ and $s=1,\dots,p'-1$ label highest weight representations.

For simplicity we only deal with the ``vacuum'' character, determined
by $|p'r-ps|=1$. The following polynomial analogues of the vacuum Virasoro
characters have arisen in the context of statistical 
mechanics~\cite{ABF84,FB85} and partition theory~\cite{ABBBFV87},
\begin{equation}\label{B}
B_L(p,p';q)=
\sum_{j=-\infty}^{\infty}\Bigl\{
q^{j(pp'j+1)}\qbin{2L}{L-p'j}
-q^{(pj+r)(p'j+s)}\qbin{2L}{L-p'j-s}\Bigr\}.
\end{equation}
The polynomials $B_L(p,p')$ are known to be related to the minimal conformal 
field theory $M(p,p')$ perturbed by the operator $\phi_{1,3}$.

Recently, very different, so-called fermionic representations
for the above polynomials have been obtained by Berkovich, McCoy and Schilling
using continued fractions~\cite{BM96,BMS98}.
Assume $p<p'<2p$, $\gcd(p,p')=1$ and define non-negative integers $n$ and 
$\nu_0,\dots,\nu_n$ by the continued fraction expansion
\begin{equation*}
\frac{p}{p'-p}=\nu_0+\cfrac{1}{\nu_1+\cfrac{1}{\nu_2+\ldots
+\cfrac{1}{\nu_n+2}}}=[\nu_0,\dots,\nu_{n-1},\nu_n+2].
\end{equation*}
Using $n$ and $\nu_j$,  set 
\begin{equation}\label{tm}
t_m=\sum_{j=0}^{m-1}\nu_j \qquad 1\leq m\leq n \qquad \text{and} \qquad
d=\sum_{j=0}^n \nu_j.
\end{equation}
The $t_m$ and $d$ are used to define a 
fractional incidence matrix $\I$ and a fractional
Cartan-type matrix $2B=2I-\I$
(with $I$ the $d$ by $d$ unit matrix) as follows
\begin{equation}\label{IB}
\I_{i,j} = \begin{cases}
\delta_{i,j+1} + \delta_{i,j-1} & \text{for
$1 \leq i<d,~ i \neq t_m$}, \\
\delta_{i,j+1} + \delta_{i,j} - \delta_{i,j-1} & \text{for
$i=t_m,~1\leq m\leq n-\delta_{\nu_n,0}$}, \\
\delta_{i,j+1} + \delta_{\nu_n,0} \delta_{i,j} & \text{for $i=d$}.
\end{cases}
\end{equation}
When $p'=p+1$ the incidence matrix $\I$ has components 
$\I_{i,j}=\delta_{|i-j|,1}$ ($i,j=1,\dots,p-2$), so that $2B$ corresponds 
to the Cartan matrix of the Lie algebra A$_{p-3}$.
When $p=2k-1$ and $p'=2k+1$ the matrix $\I$
has components $\I_{i,j}=\delta_{|i-j|,1}+
\delta_{i,j}\delta_{i,k-1}$ ($i,j=1,\dots,k-1$), so that $2B$ 
corresponds to the Cartan-type matrix of the tadpole graph of $k-1$ nodes.

Using the above definition, the fermionic representation for the bounded 
Virasoro characters with $p<p'<2p$ can be given as
\begin{equation}\label{F}
F_L(p,p';q)=\sum_{m\in 2\Integer^d}
q^{\frac{1}{2}mBm}\prod_{j=1}^d\qbin{L\delta_{j,1}+\frac{1}{2}(\I m)_j}{m_j}.
\end{equation}
Here we use the notations $vMw=\sum_{j,k} v_j M_{j,k}w_k$,
$(M v)_j=\sum_k M_{j,k}v_k$ and $(vM)_j=\sum_k v_k M_{k,j}$.
These conventions are important since, generally, 
$M\,(=\I,B)$ is not a symmetric matrix.
The general form \eqref{F} for $F_L(p,p')$ can be found in 
refs.~\cite{BM96,BMS98} (see also \cite{FLW97}). The important special
cases $(p,p')=(p,p+1)$ and $(2k-1,2k+1)$ were proven prior to this is in 
refs.~\cite{Berkovich94, Warnaar96a} and ref.~\cite{FQ95}, respectively.

The expression for $F_L(p,p';q)$ with $p'>2p$ follows from the duality 
transformation
\begin{equation}\label{Fduality}
F_L(p,p';1/q)=q^{-L^2} F_L(p'-p,p';q).
\end{equation}
To obtain fermionic character formulas for $\chi_{r,s}^{(p,p')}(q)$ with
$|p'r-ps|=1$  one simply lets $L$ tend to infinity in \eqref{F}.

Before we proceed to use the identity
\begin{equation}\label{BF}
F_L(p,p';q)=B_L(p,p';q)
\end{equation}
to derive trinomial identities let us comment on the convention of writing
$2B$ for a Cartan-type matrix in the above formulas.
This has its origin in the work of ref.~\cite{DKKMM93} where, in more 
general situations, the matrix $B$ has a (non-trivial) tensor product 
structure, $B=b_1\otimes b_2$. In the identities of this section the
matrix $b_1$ is simply the inverse of the A$_1$ Cartan matrix, 
$(b_1)=(\frac{1}{2})$. In section~\ref{secTA}, however, we indeed encounter 
a different situation, $b_1$ being the (still trivial) Cartan-type matrix 
of the tadpole graph with just a single node, so that $b_1=(1)$.

\subsection{$q$-Trinomial identities I}
We start we the $q$-binomial identity \eqref{BF} for $p<p'<2p$
assuming that $d\geq 2$. Applying equation \eqref{TB}, 
with $C_{L,k}$ given by \eqref{C}, we find
\begin{eqnarray*}
\lefteqn{\sum_{j=-\infty}^{\infty}\Bigl\{
q^{\frac{1}{2}j(p'(2p-p')j+2)}T_0(L,p'j)
-q^{\frac{1}{2}((2p-p')j+2r-s)(p'j+s)}T_0(L,p'j+s)\Bigr\}}\notag \\
\ng\ng &=&\ng\ng \sum_{k=0}^L(-1)^{L-k} q^{\binom{L-k}{2}-\frac{1}{2}L^2}
\qbins{L}{k} F_k(p,p';q) \notag\\
\ng\ng &=&\ng\ng \sum_{m\in 2\Integer^d}
q^{\frac{1}{2}(mBm-L^2)}\Bigl(\prod_{j=2}^d
\qbins{\frac{1}{2}(\I m)_j}{m_j}\Bigr)
\sum_{k=0}^L(-1)^k q^{\binom{k}{2}}
\qbins{L}{k}\qbins{L-k+\frac{1}{2}(\I m)_1}{m_1} \notag\\
\ng\ng &\stackrel{\text{by }\eqref{qcv2}}{=}&\ng\ng
q^{\frac{1}{2}L^2}\sum_{m\in 2\Integer^d}
q^{\frac{1}{2}mBm-L(Bm)_1}
\qbins{\frac{1}{2}(\I m)_1}{m_1-L}
\prod_{j=2}^d
\qbins{\frac{1}{2}(\I m)_j}{m_j} \notag\\
\ng\ng &=&\ng\ng q^{\frac{1}{4}L^2\I_{1,1}}
\sum_{m+Le_1\in 2\Integer^d}
q^{\frac{1}{2}mBm+\frac{1}{2}L(mB-Bm)_1}
\prod_{j=1}^d
\qbins{\frac{1}{2}L\I_{j,1}+\frac{1}{2}(\I m)_j}{m_j},
\end{eqnarray*}
with $e_j$, ($j=1,\dots,d$) the standard unit vectors 
in $\Integer^d$.
We now have to distinguish two cases according to whether
$\nu_0=1$ (so that $3p/2<p'<2p$) or $\nu_0>1$ (so that $p<p'\leq 3p/2$).
In the latter case $\I_{1,j}=\I_{j,1}=\delta_{1,j-1}$,
and we obtain the following polynomial identities.
\begin{proposition}\label{prop1}
For integers $p,p'$ with $p<p'\leq 3p/2$ and $\gcd(p,p')=1$ let
integers $1\leq r<p$ and $1\leq s<p'$ be fixed by $|p'r-ps|=1$ and
let $\I$ and $B$ be defined by \eqref{tm} and \eqref{IB}.
Then the following polynomial identity holds for $L\in\Integer$,
\begin{multline*}
\sum_{m+Le_1\in 2\Integer^d}q^{\frac{1}{2}mBm}\prod_{j=1}^d
\qbin{\frac{1}{2}L\delta_{j,2}+\frac{1}{2}(\I m)_j}{m_j} \\
=\sum_{j=-\infty}^{\infty}\Bigl\{
q^{\frac{1}{2}j(p'(2p-p')j+2)}T_0(L,p'j)
-q^{\frac{1}{2}((2p-p')j+2r-s)(p'j+s)}T_0(L,p'j+s)\Bigr\}.
\end{multline*}
\end{proposition}
The admissible pairs $(p,p')=(3,4)$ and $(p,p')=(2,3)$ have been
neglected in our derivation due to the assumption that $d\geq 2$.
These two cases can be treated in similar fashion and when
$(p,p')=(3,4)$ the left-hand side is $1$ for $L$ even and $0$ for $L$ odd.
When $(p,p')=(2,3)$, in which case $F_L(2,3;q)=1$, the left-hand side
becomes $\delta_{L,0}$.
All of the identities of proposition~\ref{prop1} have been derived before,
and for $p'=p+1$ they were first found by 
Schilling~\cite{Schilling96a,Schilling96b}. The more general case can be
found in ref.~\cite{SW97}.

Next we treat the case $\nu_0=1$. When this occurs
$\I_{1,j}=\delta_{j,1}-\delta_{1,j-1}$ and
$\I_{j,1}=\delta_{j,1}+\delta_{1,j-1}$
and we obtain the following polynomial identities.
\begin{proposition}\label{prop2}
For integers $p,p'$ with $3p/2<p'<2p$ and $\gcd(p,p')=1$ let
integers $1\leq r<p$ and $1\leq s<p'$ be fixed by $|p'r-ps|=1$
and let $\I$ and $B$ defined by \eqref{tm} and \eqref{IB}.
Then the following polynomial identity holds for $L\in\Integer$,
\begin{multline*}
\sum_{m+Le_1\in 2\Integer^d}
q^{\frac{1}{4}L(L-2m_2)+\frac{1}{2}mBm}\prod_{j=1}^d
\qbin{\frac{1}{2}L(\delta_{j,1}+\delta_{j,2})+\frac{1}{2}(\I m)_j}{m_j} \\
=\sum_{j=-\infty}^{\infty}\Bigl\{
q^{\frac{1}{2}j(p'(2p-p')j+2)}T_0(L,p'j)
-q^{\frac{1}{2}((2p-p')j+2r-s)(p'j+s)}T_0(L,p'j+s)\Bigr\}.
\end{multline*}
\end{proposition}
The case $(p,p')=(3,5)$ has again escaped a proper derivation, but has in
fact been treated previously, corresponding to identity \eqref{RR2} 
with $q$ replaced by $1/q$.
Apart from this special case due to Andrews~\cite{Andrews90a},
the identities of proposition~\ref{prop2}
have been proved by Berkovich, McCoy and Orrick~\cite{BMO96,BM97} for
$(p,p')=(2\nu+1,4\nu)$ and were conjectured for general $p$ and $p'$
by Berkovich, McCoy and Pearce (equation (8.8) of \cite{BMP98}).

\subsection{$q$-Trinomial identities II}
Our starting point for deriving $q$-trinomial identities is again
equation \eqref{BF} but this time we rely on \eqref{tBe}. This
implies that \eqref{BF} with $L$ replaced by $k$, multiplied
by $q^{k^2}\qbin{L}{2k}$ and summed over $k$ yields
\begin{multline}\label{even}
\sum_{k\geq 0}q^{k^2}\qbin{L}{2k}F_k(p,p';q)\\
=\sum_{j=-\infty}^{\infty}\Bigl\{
q^{j(p'(p+p')j+1)}\qbin{L}{2p'j}_2
-q^{(p'j+s)((p+p')j+r+s)}\qbin{L}{2p'j+2s}_2\Bigr\}.
\end{multline}
To transform this into explicit polynomial identities we need to
distinguish between $p<p'<2p$ and $p'>2p$.

First assume that $p<p'<2p$.
After substituting expression \eqref{F} for $F_L$, the
left side of \eqref{even} is
\begin{equation*}
\sum_{k\geq 0} \sum_{m\in 2\Integer^d}
q^{k^2+\frac{1}{4}mBm}\qbin{L}{2k}
\prod_{j=1}^d\qbin{k\delta_{j,1}+\frac{1}{2}(\I m)_j}{m_j}.
\end{equation*}
By the $q$-Chu--Vandermonde summation~\eqref{qcv1},
with $L\to L-k+m_1/2$, $a\to k-m_1/2$ and $b\to -k-m_1/2$,
this can be rewritten as
\begin{equation*}
\sum_{i,k\geq 0}\sum_{m\in 2\Integer^d}
q^{i(i-k-\frac{1}{2}m_1)+k^2+\frac{1}{4}mBm}
\qbins{L-k+\frac{1}{2}m_1}{i}\qbins{k-\frac{1}{2}m_1}{2k-i}
\prod_{j=1}^d\qbins{k\delta_{j,1}+\frac{1}{2}(\I m)_j}{m_j}.
\end{equation*}
Replacing $m_j\to m_{j+2}$, $j=1,\dots,d$, 
followed by $k\to (m_1+m_2)/2$ and $i\to m_1$ yields
\begin{multline*}
\sideset{}{'}\sum
q^{\frac{1}{4}(3m_1^2+m_2^2-2m_1 m_3)
+\frac{1}{4}\sum_{j,k=1}^d m_{j+2}B_{j,k}m_{k+2}} \\
\times 
\qbins{L-\frac{1}{2}(m_1+m_2-m_3}{m_1}\qbins{\frac{1}{2}(m_1+m_2-m_3)}{m_2}
\prod_{j=1}^d
\qbins{\frac{1}{2}(m_1+m_2)\delta_{j,1}+\frac{1}{2}\sum_{k=1}^d
\I_{j,k}m_{k+2}}{m_{j+2}},
\end{multline*}
where the primed sum denotes a sum over $m\in\Integer^{d+2}$ such
that $m_1+m_2$ and $m_3\dots,m_{d+2}$ are all even.

Now define a new incidence matrix $\I'$ and Cartan-type matrix 
$2B'=2I-\I'$ of dimension $d'=d+1$ by replacing the continued 
fraction expansion $[\nu_0,\dots,\nu_n+2]$ by $[1,\nu_0,\dots,\nu_n+2]$,
so that $\I'$ becomes the incidence matrix corresponding to the continued 
fraction expansion of $p'/p$. Also define $\I''$ and $2B''=2I-\I''$ of 
dimension $d''=d+2$ as
\begin{equation}\label{Bpp}
\I_{i,j}''=\begin{cases}
-\delta_{i,1}\delta_{j,1}+\delta_{i,2}+\delta_{i,3}-\delta_{j,2}
+\delta_{j,3} & \text{for $i=1$ or $j=1$}\\
\I_{i-1,j-1}' & \text{for $i,j=2,\dots,d+2$.}
\end{cases}
\end{equation}
Then the above sequence of transformations implies the following proposition.
\begin{proposition}\label{prop3}
For integers $p,p'$ with $p<p'<2$ and $\gcd(p,p')=1$ let
integers $1\leq r<p$ and $1\leq s<p'$ be fixed by $|p'r-ps|=1$ and
let $\I''$ and $B''$ be defined by \eqref{Bpp}.
Then the following polynomial identity holds for $L\in\Integer$,
\begin{multline*}
\sideset{}{'}\sum q^{\frac{1}{4}m B''m}
\prod_{j=1}^{d''}\qbin{L\delta_{j,1}+\frac{1}{2}(\I''m)_j}{m_j} \\
=\sum_{j=-\infty}^{\infty}\Bigl\{
q^{j(p'(p+p')j+1)}\qbin{L}{2p'j}_2
-q^{(p'j+s)((p+p')j+r+s)}\qbin{L}{2p'j+2s}_2\Bigr\}.
\end{multline*}
\end{proposition}
The identities of proposition~\ref{prop3} are the $n=0$ case
of the conjectured equation (8.11) (which contains the $n=0$ instances of 
(6.19) and (8.3)) of~\cite{BMP98}, and are related to the $\phi_{2,1}$
perturbation of the minimal conformal field theory $M(p',p+p')$.

When $p'>2p$ we replace $p\to p'-p$ in \eqref{even} and use the duality 
property~\eqref{Fduality}. Hence
\begin{multline}\label{even2}
\sum_{k\geq 0}q^{2k^2}\qbin{L}{2k}F_k(p,p';q^{-1})\\
=\sum_{j=-\infty}^{\infty}\Bigl\{
q^{j(p'(2p'-p)j+1)}\qbin{L}{2p'j}_2
-q^{(p'j+s)((2p'-p)j+r+s)}\qbin{L}{2p'j+2s}_2\Bigr\}.
\end{multline}
Observe that the transformation carried out above implies
$p<p'<2p$ and $|p'(r-s)+ps|=1$.

Substituting expression \eqref{F} for $F_L$ and using \eqref{qbindual}
the left side of \eqref{even2} yields
\begin{equation*}
\sum_{k\geq 0} \sum_{m\in 2\Integer^d}
q^{k(2k-m_1)+\frac{1}{4}mBm}
\qbin{L}{2k}\prod_{j=1}^d\qbin{k\delta_{j,1}+\frac{1}{2}(\I m)_j}{m_j}.
\end{equation*}
By the $q$-Chu--Vandermonde summation~\eqref{qcv1},
with $L\to L-m_1/2$, $a\to m_1/2$, $b\to m_1/2-2k$
this can be rewritten as
\begin{equation*}
\sum_{i,k\geq0}\sum_{m\in 2\Integer^d}
q^{i(i-2k+\frac{1}{2}m_1)+k(2k-m_1)+\frac{1}{4}mBm}
\qbins{L-\frac{1}{2}m_1}{i}\qbins{\frac{1}{2}m_1}{2k-i}
\prod_{j=1}^d\qbins{k\delta_{j,1}+\frac{1}{2}(\I m)_j}{m_j}.
\end{equation*}
Replacing $m_j\to m_{j+2}$, $j=1,\dots,d$, 
followed by $t\to m_1+m_2$ and $i\to m_1$ gives
\begin{multline*}
\sideset{}{'}\sum
q^{\frac{1}{2}(m_1^2+m_2^2-m_2 m_3)+
\frac{1}{4}\sum_{j,k=1}^d m_{j+2}B_{j,k}m_{k+2}} \\
\times 
\qbins{L-\frac{1}{2}m_3}{m_1}\qbins{\frac{1}{2}m_3}{m_2}
\prod_{j=1}^d
\qbins{\frac{1}{2}(m_1+m_2)\delta_{j,1}+\frac{1}{2}\sum_{k=1}^d
\I_{j,k}m_{k+2}}{m_{j+2}},
\end{multline*}
where the primed sum again
denotes a sum over $m\in\Integer^{d+2}$ such
that $m_1+m_2$ and $m_3\dots,m_{d+2}$ are all even.

Now define a new incidence matrix $\I$ and Cartan-type matrix
$2B'=2I-\I$ of dimension $d'=d+1$ by replacing the continued
fraction expansion $[\nu_0,\dots,\nu_n+2]$ by $[\nu_0+1,\nu_1,\dots,\nu_n+2]$,
so that $\I'$ becomes the incidence matrix corresponding to the continued
fraction expansion of $p'/(p'-p)$. Also define $\I''$ and $2B''=2I-\I''$ of 
dimension $d''=d+2$ as
\begin{equation}\label{bpp}
\I_{i,j}''=\begin{cases}
\delta_{i,3}-\delta_{j,3} & \text{for $i=1$ or $j=1$}\\
\I_{i-1,j-1}' & \text{for $i,j=2,\dots,d+2$.}
\end{cases}
\end{equation}
Then the above sequence of transformations implies the following proposition.
\begin{proposition}\label{prop4}
For integers $p,p'$ with $p<p'<3p/2$ and $\gcd(p,p')=1$ let
integers $1\leq r<p$ and $1\leq s<p'$ be fixed by $|p'(r-s)r+ps|=1$ and
let $\I$ and $B$ be defined by \eqref{bpp}.
Then the following polynomial identity holds for $L\in\Integer$,
\begin{multline*}
\sideset{}{'}\sum q^{\frac{1}{4}m B''m}
\prod_{j=1}^{d''}\qbin{L\delta_{j,1}+\frac{1}{2}(\I''m)_j}{m_j} \\
=\sum_{j=-\infty}^{\infty}\Bigl\{
q^{j(p'(2p'-p)j+1)}\qbin{L}{2p'j}_2
-q^{(p'j+s)((2p'-p)j+r+s)}\qbin{L}{2p'j+2s}_2\Bigr\}.
\end{multline*}
\end{proposition}
The identities of proposition~\ref{prop4}, which are related to
the $\phi_{2,1}$ perturbation of the conformal field theory $M(p',2p'-p)$,
were conjectured in ref.~\cite{BMP98} (as equation (6.9)). 
For $p=p'-1$ a proof using recurrences was recently given in~\cite{BM98}.

\subsection{$q$-Trinomial identities III}\label{secTA}
There are of course many more $q$-trinomial identities that can be derived
using the techniques of the previous sections. Our final application
is to show that in some cases a bit more ingenuity is required to arrive
at the desired result.
The identities we set out to prove here were again conjectured by
Berkovich, McCoy and Pearce (equation (9.4) of~\cite{BMP98})
and are interesting as they contain the (polynomial) Rogers--Ramanujan 
identity~\eqref{RR2} as simplest case.
It also provides an example for which the matrix $B=b_1\otimes b_2$ (in the
proposition below denoted $C_n$) 
of section~\ref{secBM} has $b_1=(1)$ and not $(\frac{1}{2})$.
\begin{proposition}\label{prop5}
For $n\geq 1$, let $C_n$ be the Cartan matrix of A$_n$.
Then for all $L\in\Integer$
\begin{multline}\label{TA}
\sum_{m\in\Integer^n} 
q^{\frac{1}{2}mC_n m} \prod_{j=1}^n \qbin{L\delta_{j,1}+m_j-(C_n m)_j}{m_j} \\
=\sum_{j=-\infty}^{\infty}\Bigl\{
q^{\frac{1}{2}j((n+3)(n+4)j+2)}\qbin{L}{(n+4)j}_2
-q^{\frac{1}{2}((n+3)j+2)((n+4)j+2)}\qbin{L}{(n+4)j+2}_2\Bigr\}.
\end{multline}
\end{proposition}
Letting $L$ tend to infinity using \eqref{trinlim} and \eqref{qbinlim}
this yields the following Virasoro-character identities.
\begin{corollary}
For $n\geq 1$ and $|q|<1$,
\begin{equation}\label{TAinf}
\sum_{m\in\Integer^n} 
\frac{q^{\frac{1}{2}mC_n m}}{(q)_{m_1}} \prod_{j=2}^n \qbin{m_j-(C_n m)_j}{m_j}
=\begin{cases} \chi^{((n+3)/2,n+4)}_{1,2}(q) & \text{$n$ odd} \\[1.5mm]
\chi^{((n+4)/2,n+3)}_{1,2}(q) & \text{$n$ even.}
\end{cases}
\end{equation}
\end{corollary}
In ref.~\cite{BMP98} the identities \eqref{TA} and \eqref{TAinf} were 
associated with the $\phi_{2,1}$ perturbation of the conformal field theories
$M((n+4)/2,n+3)$ when $n$ is odd and the $\phi_{1,5}$ perturbation of 
$M((n+3)/2,n+4)$ when $n$ is even.
\begin{proof}
The corollary betrays a hidden parity dependence of \eqref{TA}
which also plays a role in the proof.
Treating $n$ being odd first we set $n=2k-1$. The left-hand side
of \eqref{TA} then reads
\begin{equation}\label{eq2k1}
\sum_{m\in\Integer^{2k-1}}
q^{\frac{1}{2}m C_{2k-1} m}\prod_{j=1}^{2k-1}
\qbin{\frac{1}{2}L\delta_{j,1}+m_{j-1}-m_j+m_{j+1}}{m_j},
\end{equation}
with the convention that $m_0=L/2$ and $m_{2k}=0$.
We eliminate the variables $m_{2j-1}$, $j=1,\dots,k$ in favour of new 
variables $M_1,\dots,M_k$ defined as
\begin{equation*}
m_{2j-1}=m_{2j-2}-\frac{1}{2}(M_j-M_{j+1}),
\end{equation*}
where $M_{k+1}=0$.
If after this replacement we relabel $m_{2j}$ to $m_j$ 
for $j=1,\dots,k$ (so that $m_k=0$), expression \eqref{eq2k1} becomes
\begin{multline}\label{mM}
\sum_{M+Le_1\in 2\Integer^k}
q^{\frac{1}{4}(L(L-2M_1)+M_1^2+
\sum_{i,j=2}^k M_i (C_{k-1})_{i,j} M_j)} \\
\times \sum_{m_1,\dots,m_{k-1}} 
q^{\sum_{j=1}^{k-1}(M_{j+1}-m_j)(m_{j-1}-m_j-\frac{1}{2}(M_j-M_{j+2}))} \\
\times\qbins{m_0+m_1+\frac{1}{2}(M_1-M_2)}{m_0-\frac{1}{2}(M_1-M_2)}
\prod_{j=1}^{k-1}
\qbins{m_{j-1}-\frac{1}{2}(M_j-M_{j+2})}{m_j}
\qbins{m_{j+1}+\frac{1}{2}(M_{j+1}-M_{j+2})}{m_j-\frac{1}{2}(M_{j+1}-M_{j+2})}.
\end{multline}
This allows for successive summation over $m_{k-1},\dots,m_1$ by 
the $q$-Saalsch\"utz sum \eqref{qS}.  When summing over $m_j$ we take 
\eqref{qS} with $L\to m_{j-1}-(M_j-M_{j+2})/2$,
$a\to (M_{j+1}-M_{j+2})/2$, $b\to -(M_{j+1}+M_{j+2})/2$,
$c\to (M_j-M_{j+1})/2$ (for $j\geq 2$) and
$c\to m_0+(M_1-M_2)/2$ (for $j=1$).
As a result \eqref{mM} collapses into
\begin{equation}\label{Fred}
\sum_{M+Le_1\in 2\Integer^k}
q^{\frac{1}{4}L(L-2M_1)+\frac{1}{2}M B M}\prod_{j=1}^k 
\qbin{\frac{1}{2}L(\delta_{j,1}+\delta_{j,2})+\frac{1}{2}(\I M)_j}{M_j},
\end{equation}
with matrices $\I$ and $2B=2I-\I$ defined in equations~\eqref{tm} 
and \eqref{IB} corresponding to the continued fraction expansion
$(k+2)/(k+1)=[1,k-1]$, i.e.,
\begin{equation*}
\I_{i,j}=\begin{cases}
\delta_{i,1}\delta_{j,1}+\delta_{i,2}-\delta_{j,2} & 
\text{for $i=1$ or $j=1$} \\
\delta_{i,j-1}+\delta_{i,j+1} & \text{for $i,j=2,\dots,k$.}
\end{cases}
\end{equation*}
The last part of the proof consist of the observation that the identity
obtained by equating \eqref{Fred} with the right-hand side of \eqref{TA} 
(with $n=2k-1$) is nothing but the identity of proposition~\ref{prop2} with 
$(p,p')=(k+2,2k+3)$ after letting $q\to 1/q$. This is readily seen using
\eqref{qtT} and \eqref{qbindual}.

Next we deal with $n$ being even setting $n=2k$.
The left-hand side of \eqref{TA} then is 
\begin{equation}\label{eq2k}
\sum_{m\in\Integer^{2k}}
q^{\frac{1}{2}m C_{2k} m}\prod_{j=1}^{2k}
\qbin{L\delta_{j,1}+m_{j-1}-m_j+m_{j+1}}{m_j},
\end{equation}
where $m_0=m_{2k+1}=0$.
We eliminate the variables $m_{2j}$, $j=1,\dots,k$, introducing new
variables $M_0,\dots,M_{k-1}$ by
\begin{equation*}
m_{2j}=m_{2j-1}-\frac{1}{2}(M_{j-1}-M_j),
\end{equation*}
where $M_k=0$. After this replacement we shift
$m_{2j-1}\to m_j$ for $j=1,\dots,k$
so that expression \eqref{eq2k} becomes
\begin{multline*}
\sum_{\substack{M_0,\dots,M_{k-1}\\M_j\text{ even}}}
q^{\frac{1}{4}(M_0^2+\sum_{i,j=1}^{k-1}M_i (C_{k-1})_{i,j}M_j)} \\
\times \sum_{m_1,\dots,m_k} q^{m_1(m_1-\frac{1}{2}(M_0+M_1))+
\sum_{j=2}^k 
(M_{j-1}-m_j)(m_{j-1}-m_j-\frac{1}{2}(M_{j-2}-M_j))} \\
\times\qbins{L-\frac{1}{2}(M_0-M_1)}{m_1}
\Bigl(\prod_{j=2}^k
\qbins{m_{j-1}-\frac{1}{2}(M_{j-2}-M_j)}{m_j}\Bigr)
\prod_{j=1}^k
\qbins{m_{j+1}+\frac{1}{2}(M_{j-1}-M_j)}
{m_j-\frac{1}{2}(M_{j-1}-M_j)}.
\end{multline*}
We now sum over $m_k,\dots,m_3$ by successive application of the 
$q$-Saalsch\"utz sum \eqref{qS}.
When summing over $m_j$ we take \eqref{qS}
with $L\to m_{j-1}-(M_{j-2}-M_j)/2$, $a\to (M_{j-1}-M_j)/2$,
$b\to -(M_{j-1}+M_j)/2$ and $c\to (M_{j-2}-M_{j-1})/2$.
The final sum over $m_1$ follows from \eqref{qcv1}
with $L\to L-(M_0-M_1)/2$, $a\to (M_0-M_1)/2$ and $b\to -(M_0+M_1)/2$.
Setting $M_0\to 2i$ the resulting expression is 
\begin{equation*}
\sum_{i\geq 0} q^{i^2}\qbin{L}{2i}\sum_{M\in 2\Integer^{k-1}}
q^{\frac{1}{4}MC_{k-1} M}
\prod_{j=1}^{k-1}\qbin{i\delta_{j,1}+M_j-\frac{1}{2}(C_{k-1}M)_j}{M_j}.
\end{equation*}
Equating this with the right-hand side of \eqref{TA} for $n=2k$ we 
recognize identity \eqref{even} with $(p,p')=(k+1,k+2)$.
\end{proof}

\section{The trinomial Bailey lemma}\label{secBl}
In this final section of our paper we formulate some of our
results in the language of Bailey pairs. As we will see, the
connection coefficients obtained in section~\ref{seccon} provide a
very elementary proof of the trinomial analogue of Bailey's
lemma recently obtained by Andrews and Berkovich~\cite{AB98}.

First some definitions are needed. In subsequent formulas
$T_n(L,a)/(q)_L$ will be abbreviated to $Q_n(L,a)$.
\begin{definition}
A pair of sequences $\alpha=\{\alpha_L\}_{L\geq 0}$ and
$\beta=\{\beta_L\}_{L\geq 0}$ that satisfies
\begin{equation*}
\beta_L = \sum_{r=0}^L \frac{\alpha_r}{(q)_{L-r}(aq)_{L+r}}
\end{equation*}
forms a (binomial) Bailey pair relative to $a$.
\end{definition}
\begin{definition}
A pair of sequences $A=\{A_L\}_{L\geq 0}$ and
$B=\{B_L\}_{L\geq 0}$ that satisfies
\begin{equation*}
B_L=\sum_{r=0}^L Q_n(L,r)A_r
\end{equation*}
forms a trinomial Bailey pair relative to $n$.
\end{definition}
The Bailey lemma~\cite{Bailey49} and trinomial Bailey lemma~\cite{AB98}
can now be stated as the following summation formulas.
\begin{lemma}\label{BaileyLemma}
Let $(\alpha,\beta)$ be a Bailey pair relative to $a$. Then
\begin{multline*}
\sum_{L=0}^M
\frac{(\rho_1)_L (\rho_2)_L (aq/\rho_1 \rho_2)^L \alpha_L}
{(aq/\rho_1)_L (aq/\rho_2)_L(q)_{M-L}(aq)_{M+L}} \\
=\sum_{L=0}^M
\frac{(\rho_1)_L (\rho_2)_L (aq/\rho_1 \rho_2)^L
(aq/\rho_1 \rho_2)_{M-L} \beta_L}
{(aq/\rho_1)_M (aq/\rho_2)_M (q)_{M-L}}.
\end{multline*}
\end{lemma}
\begin{lemma}\label{trinBaileyLemma}
Let $(A,B)$ form a trinomial Bailey pair relative to $0$. Then
\begin{equation}\label{tbl1}
\sum_{L=0}^M (-1)_L \, q^{\frac{1}{2}L} B_L =
(-1)_{M+1} \sum_{L=0}^M q^{\frac{1}{2}L}A_L\frac{Q_1(M,L)}{1+q^L}\, .
\end{equation}
If $(A,B)$ is a trinomial Bailey pair relative to $1$, then
\begin{multline*}
\sum_{L=0}^M \big(-q^{-1}\big)_L\, q^L B_L  \\
= (-1)_M \sum_{L=0}^M A_L \biggl\{ Q_1(M,L)
-\frac{Q_1(M-1,L+1)}{1+q^{-L-1}}
-\frac{Q_1(M-1,L-1)}{1+q^{L-1}} \biggr\}.
\end{multline*}
\end{lemma}

Before we translate the results of section~\ref{seccon} in the 
language of Bailey pairs let us point out that the connection coefficients
between $q$-binomials and $q$-trinomials can be applied to yield a very
simple proof of the trinomial Bailey lemma.
At the heart of the proof of lemma~\ref{trinBaileyLemma} is the following
identity derived in \cite{AB98} by a considerable amount of work
\begin{equation}\label{Tbl}
T_0(L,a)=q^{\frac{1}{2}(a-L)}\Bigl\{
\frac{1+q^L}{1+q^a} \, T_1(L,a)-
\frac{1-q^L}{1+q^a} \, T_1(L-1,a)\Bigr\}.
\end{equation}
To see, for example, that this implies \eqref{tbl1} we multiply \eqref{Tbl}
by $q^{L/2}(-1)_L/(q)_L$ and sum over $L$ from $a$ to $M$. On the right-hand
side all but one term cancels, so that
\begin{equation*}
\sum_{L=a}^M q^{\frac{1}{2}L}(-1)_L Q_0(L,a)
=\frac{q^{\frac{1}{2}a}}{1+q^a} (-1)_{M+1}Q_1(M,a),
\end{equation*}
which obviously implies \eqref{tbl1}.

By equations~\eqref{TB}--\eqref{Cp} equation~\eqref{Tbl} is proved if we can
show its validity when multiplied by $C_{M,L}(a)$ and summed over $L$. Doing
this and using \eqref{BT}, one finds (replacing $L\to k$ and $M\to L$)
\begin{align*}
\qbin{2L}{L-a}&=\sum_{k=0}^L q^{\binom{k}{2}-\binom{a}{2}}\qbin{L}{k}
\Bigl\{
\frac{1+q^k}{1+q^a} \, T_1(k,a)-
\frac{1-q^k}{1+q^a} \, T_1(k-1,a)\Bigr\} \\
&=\sum_{k=0}^L q^{\binom{k}{2}-\binom{a}{2}}
\Bigl\{
\frac{1+q^k}{1+q^a}\qbin{L}{k}-
q^k \frac{1-q^{k+1}}{1+q^a}\qbin{L}{k+1}\Bigr\} T_1(k,a)\\
&=\sum_{k=0}^L q^{\binom{k}{2}-\binom{a}{2}}
\frac{1+q^L}{1+q^a}\qbin{L}{k} T_1(k,a)
\end{align*}
But the extremes of this string of equations is nothing but equation 
\eqref{BT1} with $D'_{L,k}(a)$ given by equation \eqref{Dp} of 
lemma~\ref{lemma2}, establishing \eqref{Tbl}.

We now give a series of lemmas which are all straightforward
consequences of the results of section~\ref{seccon}.
\begin{lemma}\label{lemAB}
Let $(\alpha,\beta)$ be a Bailey pair relative to $1$. Then
\begin{align*}
A_L &= q^{-\frac{1}{2}L^2}\alpha_L \notag \\
B_L &=\sum_{k=0}^L\frac{(-1)^{L-k}q^{\binom{L-k}{2}-\frac{1}{2}L^2}(q)_{2k}}
{(q)_k(q)_{L-k}} \,\beta_k
\end{align*}
is a trinomial Bailey pair relative to $0$ and
\begin{align*}
A_L &= \frac{q^{-\binom{L}{2}}}{1+q^L}\, \alpha_L\notag \\
B_L &=\sum_{k=0}^L\frac{(-1)^{L-k}q^{\binom{L-k}{2}-\binom{L}{2}}(q)_{2k}}
{(1+q^k)(q)_k(q)_{L-k}} \,\beta_k
\end{align*}
is a trinomial Bailey pair relative to $1$.
\end{lemma} 
The converse statement is as follows.
\begin{lemma}
Let $(A(n),B(n))$ be a trinomial Bailey pair relative to $n$. Then,
\begin{align*}
\alpha_L &= q^{\frac{1}{2}L^2}A_L(0) \notag \\
\beta_L &=\frac{(q)_L}{(q)_{2L}}\sum_{k=0}^L
\frac{q^{\frac{1}{2}k^2}}{(q)_{L-k}} \,B_k(0)
\end{align*}
and
\begin{align*}
\alpha_L &= q^{\binom{L}{2}}(1+q^L)A_L(1) \notag \\
\beta_L &=\frac{(q)_L}{(q)_{2L}}(1+q^L)\sum_{k=0}^L
\frac{q^{\binom{k}{2}}}{(q)_{L-k}} \,B_k(1)
\end{align*}
are Bailey pairs relative to $1$.
\end{lemma}
Lemma~\ref{lemAB} is to be compared with the following
result of ref.~\cite{Warnaar98}.
\begin{lemma}
Let $\ell$ be a non-negative integer and $(\alpha,\beta)$ a Bailey pair 
relative to $a=q^{\ell}$. Then
\begin{align*}
A_L&=\begin{cases}
\alpha_{\frac{1}{2}(L-\ell)} &\text{for $L=\ell,\ell+2,\dots$} \\
0 &\text{otherwise} 
\end{cases}  \notag \\
B_L&=\begin{cases} \displaystyle
\sum_{k=0}^{\lfloor (L-\ell)/2\rfloor}
\frac{q^{\frac{1}{2}(L-\ell-2k)(L-\ell-2k-n)}}
{(q)_{\ell}(q)_{L-\ell-2k}} \beta_k & \text{for $L\geq \ell$} \\
0 &\text{otherwise}
\end{cases}
\end{align*}
forms a trinomial Bailey pair relative to $n$.
\end{lemma}

\appendix
\section{Some $q$-binomial formulas}
In this appendix we list some standard $q$-binomial identities that are 
repeatedly used in the main text.

The following three formulas all hold for integers $a,b,L$ such that 
$a,L\geq 0$,
\begin{equation}\label{qcv1}
\sum_{k=0}^L q^{k(k+b)}\qbin{L}{k}\qbin{a}{k+b}=\qbin{a+L}{b+L},
\end{equation}
\begin{equation}\label{qcv2}
\sum_{k=0}^L (-1)^k q^{\binom{k}{2}}\qbin{L}{k}\qbin{L+a-k}{b}=
q^{L(L+a-b)}\qbin{a}{b-L},
\end{equation}
\begin{equation}\label{qcv3}
\sum_{k=0}^L (-1)^k q^{\binom{k}{2}+k(b-L+1)}\qbin{L}{k}\qbin{L+a-k}{b}=
\qbin{a}{b-L}.
\end{equation}
The first two equation are specializations of the $q$-Chu--Vandermonde
sum (II.7) of \cite{GR90} and the last equation is a specialization of
the $q$-Chu--Vandermonde sum (II.6) of \cite{GR90}. Identity~\eqref{qcv2}
is also given in \cite{Andrews76} as equation (3.3.10).
A useful specialization of the $q$-Saalsch\"utz sum ((II.12) of \cite{GR90})
is given by
\begin{equation}\label{qS}
\sum_{k=0}^L q^{(a-b-k)(L-k)}
\qbin{L}{k}\qbin{a}{k+b}\qbin{k+c}{a+L}=
\qbin{c}{b+L}\qbin{c-b}{a-b},
\end{equation}
true for integers $a,b,c,L$ such that $a,c,L\geq 0$. This is equation
(3.3.11) of \cite{Andrews76}.
Finally we list the elementary results
\begin{equation}\label{qbinrec}
\qbin{L}{a}=\qbin{L-1}{a-1}+q^a \qbin{L-1}{a} \qquad \text{for }L,a\geq 0,
L+a\neq 0,
\end{equation}
\begin{equation}\label{qbinlim}
\lim_{L\to\infty}\qbin{L}{a}=\frac{1}{(q)_a}
\end{equation}
and
\begin{equation}\label{qbindual}
\qbin{L}{a}_{1/q}=q^{-a(L-a)}\qbin{L}{a}_q.
\end{equation}

\subsection*{Acknowledgements}
I thank Alexander Berkovich and Anne Schilling for helpful
discussions. 
This work has been supported by a fellowship of the Royal
Netherlands Academy of Arts and Sciences.

\end{document}